\documentclass[11pt,reqno]{amsart}
\usepackage{fullpage}

\renewcommand{\le}{\leqslant}
\renewcommand{\ge}{\geqslant}

\renewcommand{\setminus}{\smallsetminus}

\renewcommand{\subset}{\subseteq}

\newcommand{\n}{\{1,\ldots,n\}}

\newcommand{\e}{\varepsilon}

\newcommand{\sign}{\mathrm{sign}}
\usepackage{amsmath,amsthm,amsfonts,amssymb,verbatim}
\usepackage{paralist}
\usepackage{enumerate}
\usepackage[all]{xy}
\usepackage[mathscr]{euscript}
\usepackage{color}
\usepackage{upgreek}
\renewcommand{\gamma}{\upgamma}
\renewcommand{\1}{\mathbf 1}
\usepackage{graphicx,verbatim}
\usepackage{paralist}
\usepackage[breaklinks,pdfstartview=FitH]{hyperref}

\newcommand{\alert}[1]{\textbf{\color{red}
[[[#1]]]}\marginpar{\textbf{\color{red}**}}\typeout{ALERT:
\the\inputlineno: #1}}
\newcommand{\N}{\mathbb{N}}

\newcommand{\R}{\mathbb{R}}

\newcommand{\mommit}[1]{}
\newcommand{\namedref}[2]{\hyperref[#2]{#1~\ref*{#2}}}

\theoremstyle{plain}
\newtheorem{theorem}{Theorem}

\newtheorem{remark}[theorem]{Remark}

\theoremstyle{definition}

\newenvironment{RETHM}[2]{\trivlist \item[\hskip 10pt\hskip\labelsep{\bf
#1\hskip 5pt\relax\ref{#2}.}]\it}{\endtrivlist}
\newcommand{\rethm}[1]{\begin{RETHM}{Theorem}{#1}}
\newcommand{\repro}[1]{\begin{RETHM}{Proposition}{#1}}
\newcommand{\relem}[1]{\begin{RETHM}{Lemma}{#1}}
\newcommand{\recor}[1]{\begin{RETHM}{Corollary}{#1}}
\newcommand{\erethm}{\end{RETHM}}
\renewcommand{\epsilon}{\varepsilon}

\newcommand{\eqdef}{\stackrel{\mathrm{def}}{=}}

\title[Uniform nonextendability from nets]{Uniform nonextendability from nets}

\author{Assaf Naor}
\address{Mathematics Department\\ Princeton University\\ Fine Hall, Washington Road, Princeton, NJ 08544-1000, USA.}
\email{naor@math.princeton.edu}
\thanks{Supported by NSF grant CCF-0832795, BSF grant 2010021, the Packard Foundation and the Simons Foundation.}

\begin{document}
\maketitle


\begin{abstract}
It is shown that there exist Banach spaces $X,Y$, a $1$-net $\mathscr{N}$ of $X$ and a Lipschitz function $f:\mathscr{N}\to Y$ such that every $F:X\to Y$ that extends $f$ is not uniformly continuous.
\end{abstract}


\section{Introduction}

A metric space $(X,d_X)$ is said to embed uniformly into a metric space $(Y,d_Y)$ if there exists an injection $f:X\to Y$ such that both $f$ and $f^{-1}$ are uniformly continuous. $(X,d_X)$ is said~\cite{Gro93} to embed coarsely into $(Y,d_Y)$ if there exists $f:X\to Y$ and nondecreasing functions $\alpha,\beta:[0,\infty)\to [0,\infty)$ with $\lim_{t\to \infty}\alpha(t)=\infty$ such that $\alpha(d_X(x,y))\le d_Y(f(x),f(y))\le \beta(d_X(x,y))$ for every $x,y\in X$. While making no attempt to survey the very large literature on these topics, we only indicate here that in addition to their intrinsic geometric interest, uniform and coarse embeddings have important applications in areas ranging from functional analysis~\cite{BL00} to group theory and topology~\cite{Yu06}, and theoretical computer science~\cite{AKR15}.

In the context of embeddings of Banach spaces, the literature suggests that  uniform and coarse embeddings are closely related, despite dealing with infinitesimal and large-scale structures, respectively. Specifically, by~\cite{AMM85,JR06,Ran06} a Banach space $X$ embeds uniformly into a Hilbert space if and only if it embeds coarsely into a Hilbert space. Also, certain obstructions work equally well~\cite{Kal07,MN08,Kra14} for ruling out both uniform and coarse embeddings of Banach spaces. Despite this, it remains unknown whether or not the existence of a coarse embedding of a Banach space $X$ into a Banach space $Y$ implies that $X$ also embeds uniformly into $Y$. The analogous question with the roles of coarse and uniform embeddings interchanged is open as well. The only available negative result in this context treats uniform and coarse equivalences rather than embeddings: Kalton~\cite{Kal12} proved the existence of two Banach spaces $X,Y$ that are coarsely equivalent but not uniformly equivalent.

Recent work of Rosendal~\cite{Ros15} yields progress towards the above questions. It implies  that if $X$ and $Y$ are Banach spaces such that $X$ embeds uniformly into $Y$,  then $X$ also embeds coarsely into $\ell_p(Y)$ for every $p\ge 1$. As for the deduction of uniform embeddability from coarse embeddability, Rosendal's work~\cite{Ros15} implies that if $X$ and $Y$ are Banach spaces with the property that for every $1$-net $\mathscr{N}$ of $X$, every Lipschitz function $f:\mathscr{N}\to Y$ admits an extension $F:X\to Y$ that is uniformly continuous, then the existence of a coarse embedding of $X$ into $Y$ implies that $X$ embeds uniformly into $\ell_p(Y)$ for every $p\ge 1$. Rosendal therefore asked~\cite{Ros15} whether or not every pair of Banach spaces $X,Y$ has this (seemingly weak) extension property. Here we show that this is not the case.

\begin{theorem}\label{thm:no ext}
There exist two Banach spaces $(X,\|\cdot\|_X)$ and $(Y,\|\cdot\|_Y)$, a $1$-net $\mathscr{N}$ of $X$ and a Lipschitz function $f:\mathscr{N}\to Y$ such that every $F:X\to Y$ that extends $f$ is not uniformly continuous. Moreover, any $F:X\to Y$ that is uniformly continuous satisfies
\begin{equation}\label{eq:no bounded distance}
\sup_{x\in \mathscr{N}} \|F(x)-f(x)\|_Y=\infty.
\end{equation}
\end{theorem}

It remains an interesting open question to understand those pairs of Banach spaces $X,Y$ for which Rosendal's question has a positive answer, i.e., to prove theorems asserting the existence of a uniformly continuous extension of any $Y$-valued Lipschitz function that is defined on a $1$-net $\mathscr{N}$ of $X$. If the initial function $f:\mathscr{N}\to Y$ is assumed to be H\"older with sufficiently small exponent rather than Lipschitz  (which is a more stringent requirement since the minimum positive distance in $\mathscr{N}$ is at least $1$), then such an extension result holds true provided that $Y$ is superreflexive. This follows from deep work of Ball~\cite{Bal92} (see~\cite{MN13} for the precise statement that we need here). Indeed, if $Y$ is superreflexive then by the work of Pisier~\cite{Pis75} we know that $Y$ admits an equivalent norm whose modulus of uniform convexity has power-type $q$ for some $q\in [2,\infty)$. Therefore, by~\cite{Bal92,MN13}, $Y$ has metric Markov cotype $q$. Since any metric raised to the power $1/q$ has Markov type $q$, it follows from~\cite{Bal92,MN13} that every $1/q$-H\"older function from a subset of $X$ into $Y$ can be extended to a $1/q$-H\"older function defined on all $X$. The role of superreflexivity here is only through the finiteness of the metric Markov cotype of $Y$, so by~\cite{MN13} similar statements hold true when $Y$ is a $q$-barycentric metric space (in particular, if $Y$ is a Hadamard space then this holds true with $q=2$). These considerations, however, do not address Rosendal's question, where the initial function $f:\mathscr{N}\to Y$ is only assumed to be Lipschitz. At the same time, an inspection of the proof below reveals that in Theorem~\ref{thm:no ext} we can ensure for every $\alpha\in (0,1)$ that $f$ is $\alpha$-H\"older, so in general the uniformly continuous extension problem for H\"older functions on nets in Banach spaces has a negative answer.



\section{Proof of Theorem~\ref{thm:no ext}}

The proof below is a variant of the argument in Section~5 of~\cite{Nao01}, which itself uses an averaging idea that is inspired by Lemma~6 in~\cite{Lin64}.

For $p\ge 2$ let $M_{p}:\ell_2\to \ell_p$ be the Mazur map, i.e., for every $x\in \ell_2$ and $j\in \N$,
$$
M_p(x)_j\eqdef |x_j|^{\frac{2}{p}}\sign(x_j).
$$
Since $p\ge 2$, it is elementary to check that every $u,v\in \R$ satisfy
$$
\left||u|^{\frac{2}{p}}\sign(u)-|v|^{\frac{2}{p}}\sign(v)\right|^p\le 2^{p-2}|u-v|^2.
$$
Consequently, for every $x,y\in \ell_2$ we have.
\begin{equation}\label{eq:Mp holder}
\|M_p(x)-M_p(y)\|_p\le 2^{1-\frac{2}{p}}\|x-y\|_2^{\frac{2}{p}}\le 2\|x-y\|_2^{\frac{2}{p}}.
\end{equation}
Denote
\begin{equation}\label{eq:X,Y}
X\eqdef \Big(\bigoplus_{p=2}^\infty \ell_2\Big)_\infty\qquad \mathrm{and}\qquad Y\eqdef \Big(\bigoplus_{p=2}^\infty  \ell_p\Big)_\infty.
\end{equation}

Fix a $1$-net $\mathscr{M}$ of $\ell_2$ and denote $\mathscr{N}=\prod_{p=2}^\infty\mathscr{M}$. Then (by the definition of $X$) $\mathscr{N}$ is a $1$-net of $X$. Define $f:\mathscr{N}\to Y$ by setting for every $(x_p)_{p=2}^\infty\in X$,
$$
f((x_p)_{p=2}^\infty)\eqdef (M_p(x_p))_{p=2}^\infty.
$$
Since the minimum distance in $\mathscr{N}$ is at least $1$, it follows from~\eqref{eq:Mp holder} (and the definitions of $X$ and $Y$) that $f$ is $2$-Lipschitz.

Suppose for the purpose of obtaining a contradiction that there exists $F:X\to Y$ that is uniformly continuous and satisfies
\begin{equation}\label{eq:bounded distance}
\gamma \eqdef \sup_{x\in \mathscr{N}} \|F(x)-f(x)\|_Y<\infty.
\end{equation}
Let $\omega:[0,\infty)\to [0,\infty)$ be the modulus of uniform continuity of $F$. Thus $\omega$ is nondecreasing and $\lim_{s\to 0} \omega(s)=0$.  Write $F=(F_p)_{p=2}^\infty$, where for every integer $p\ge 2$ the mapping  $F_p:\ell_2\to \ell_p$ also has modulus of continuity that is bounded from above by $\omega$ (by the definitions of $X$ and $Y$). By~\eqref{eq:bounded distance} and the definition of $f$, for every $y\in \mathscr{M}$ we have $\|F_p(y)-M_p(y)\|_p\le \gamma$. Hence, since $\mathscr{M}$ is a $1$-net,
\begin{multline}\label{eq:close to Mp}
\sup_{x\in \ell_2} \|F_p(x)-M_p(x)\|_p\le \sup_{x\in \ell_2} \inf_{y\in \mathscr{M}}\big(\|F_p(x)-F_p(y)\|_p+\|F_p(y)-M_p(y)\|_p+\|M_p(y)-M_p(x)\|_p\big)\\\stackrel{\eqref{eq:Mp holder}}{\le}  \sup_{x\in \ell_2} \inf_{y\in \mathscr{M}}\left(\omega\left(\|x-y\|_2\right)+\gamma +2\|y-x\|_2^{\frac{2}{p}}\right)\le \omega(1)+\gamma+2.
\end{multline}

In what follows, for every $n\in \N$ we let $J_n:\ell_2^n\to \ell_2$ be the canonical embedding, i.e., $J_n(x_1,\ldots,x_n)=(x_1,\ldots,x_n,0,\ldots)$. Also, we let $Q_n:\ell_p\to \ell_p^n$ be the canonical projection, i.e., $Q_n((x_j)_{j=1}^\infty)=(x_1,\ldots,x_n)$. Given $n\in \N$, we identify a permutation $\pi\in S_n$ with its associated permutation matrix, i.e.,  $\pi x=(x_{\pi^{-1}(1)},\ldots,x_{\pi^{-1}(n)})$ for every  $x\in \R^n$. Similarly, we identify $\e\in \{-1,1\}^n$ with the corresponding diagonal matrix, i.e.,  $\e x=(\e_1x_{1},\ldots,\e_n x_{n})$ for every $x\in \R^n$.

Fix two integers $p,n\in \N$ with $p\ge 2$. Define $G_p^n:\ell_2^n\to \ell_p^n$ by
\begin{equation}\label{eq:def Gp}
\forall\, x\in \ell_2^n,\qquad G_p^n(x)\eqdef  \frac{1}{2^nn!}\sum_{\pi\in S_n}\sum_{\e\in \{-1,1\}^n}(\e\pi)^{-1} Q_n\circ F_p\circ J_n(\e\pi x).
\end{equation}
Then, because $\{\e\pi:\ (\e,\pi)\in \{-1,1\}^n\times S_n\}$ forms a group of linear operators on $\R^n$, we have $G_p^n(\e\pi x)=\e\pi G_p^n(x)$ for every $(\e,\pi)\in \{-1,1\}^n\times  S_n$ and $x\in \ell_2^n$. Since for every $A\subset \n$ and $t\in \R$ we have $\pi (t\1_A)=t\1_A$ whenever $\pi \in S_n$ fixes $A$, it follows that there exist $\alpha(t,A),\beta(t,A)\in \R$ such that $G_p^n(t\1_A)=\alpha(t,A)\1_A+\beta(t,A)\1_{\n\setminus A}$. Since $(\1_A-\1_{\n\setminus A})(t\1_A)=t\1_A$, it follows that $\beta(t,A)=-\beta(t,A)$, thus $G_p^n(t\1_A)=\alpha(t,A)\1_A$. Finally, since for every $A,B\subset \n$ of the same cardinality there exists $\pi\in S_n$ with $\pi(t\1_A)=t\1_B$, we conclude that $\alpha(t,A)$ depends only on the cardinality of $A$. In other words, there exists a sequence $\{\alpha_k(t)\}_{k=0}^n\subset \R$ such that
\begin{equation}\label{eq:indicator to indicator}
\forall\, A\subset \n,\ \forall\, t\in \R,\qquad G_p^n(t\1_A)=\alpha_{|A|}(t)\1_A.
\end{equation}

Since $Q_n\circ M_p\circ J_n(\e\pi x)=\e\pi Q_n\circ M_p\circ J_n(x)$ for every $x\in \ell_2^n$ and $(\e,\pi)\in \{-1,1\}^n\times S_n$,
\begin{align*}
\sup_{x\in \ell_2^n} \left\|G_p^n(x)-Q_n\circ M_p\circ J_n(x)\right\|_p&\stackrel{\eqref{eq:def Gp} }{=}\sup_{x\in \ell_2^n}\Big\|\frac{1}{2^nn!}\sum_{\pi\in S_n}\sum_{\e\in \{-1,1\}^n}(\e\pi)^{-1} Q_n\circ (F_p-M_p)\circ J_n(\e\pi x)\Big\|_p\\&\stackrel{\eqref{eq:close to Mp}}{\le} \omega(1)+\gamma+2.
\end{align*}
In particular, for every $k\in \{1,\ldots, n-1\}$ and $t\in (0,\infty)$ we have
\begin{multline*}\label{eq:first approx}
\left\|G_p^n\left(t\1_{\{n-k+1,\ldots,n\}}\right)-t^{\frac{2}{p}}\1_{\{n-k+1,\ldots,n\}}\right\|_p\\ =\left\|G_p^n\left(t\1_{\{n-k+1,\ldots,n\}}\right)-Q_n\circ M_p\circ J_n\circ\left(t\1_{\{n-k+1,\ldots,n\}}\right)\right\|_p\le \omega(1)+\gamma+2.
\end{multline*}
and
\begin{equation*}\label{erq:second approx}
\left\|G_p^n\left(t\1_{\{1,\ldots,k\}}\right)-t^{\frac{2}{p}}\1_{\{1,\ldots,k\}}\right\|_p=\left\|G_p^n\left(t\1_{\{1,\ldots,k\}}\right)-Q_n\circ M_p\circ J_n\circ\left(t\1_{\{1,\ldots,k\}}\right)\right\|_p\le \omega(1)+\gamma+2.
\end{equation*}
Consequently, assuming from now on that $2k\le n+1$ we have
\begin{align}\label{eq:large because of mazur}
t^{\frac{2}{p}}(2k)^{\frac{1}{p}}&=\left\|t^{\frac{2}{p}}\left(\1_{\{n-k+1,\ldots,n\}}-\1_{\{1,\ldots,k\}}\right)\right\|_p\nonumber \\\nonumber &\le \left\|G_p^n\left(t\1_{\{n-k+1,\ldots,n\}}\right)-G_p^n\left(t\1_{\{1,\ldots,k\}}\right)\right\|_p
+\left\|G_p^n\left(t\1_{\{n-k+1,\ldots,n\}}\right)-t^{\frac{2}{p}}\1_{\{n-k+1,\ldots,n\}}\right\|_p\\\nonumber &\qquad
+ \left\|G_p^n\left(t\1_{\{1,\ldots,k\}}\right)-t^{\frac{2}{p}}\1_{\{1,\ldots,k\}}\right\|_p\\&\le  \left\|G_p^n\left(t\1_{\{n-k+1,\ldots,n\}}\right)-G_p^n\left(t\1_{\{1,\ldots,k\}}\right)\right\|_p+2\omega(1)+2\gamma+4.
\end{align}

At the same time, since $G_p^n$ is obtained by averaging compositions of $Q_n\circ F_p\circ J_n:\ell_2^n\to \ell_p^n$ with isometries (of both the source space and the target space), the modulus of uniform continuity of $G_p^n$ is bounded from above by $\omega$. Hence,
\begin{multline*}
\left\|G_p^n\left(t\1_{\{n-k+1,\ldots,n\}}\right)-G_p^n\left(t\1_{\{1,\ldots,k\}}\right)\right\|_p\stackrel{\eqref{eq:indicator to indicator}}{=} \left\|\alpha_k(t)\left(\1_{\{n-k+1,\ldots,n\}}-\1_{\{1,\ldots,k\}}\right)\right\|_p= |\alpha_k(t)|(2k)^{\frac{1}{p}}\\
= k^{\frac{1}{p}}\left\|\alpha_k(t)\left(\1_{\{1,\ldots,k\}}-\1_{\{2,\ldots,k+1\}}\right)\right\|_p \stackrel{\eqref{eq:indicator to indicator}}{=}  k^{\frac{1}{p}}\left\|G_p^n\left(t\1_{\{1,\ldots,k\}}\right)-G_p^n\left(t\1_{\{2,\ldots,k+1\}}\right)\right\|_p\le k^{\frac{1}{p}} \omega(\sqrt{2} t).
\end{multline*}
In combination with~\eqref{eq:large because of mazur}, this yields the following estimate, which holds for every $p,k\in \N$ with $p\ge 2$ and $t\in (0,\infty)$.
\begin{equation}\label{eq:omega constraint}
t^{\frac{2}{p}}(2k)^{\frac{1}{p}}\le k^{\frac{1}{p}} \omega(\sqrt{2} t)+ 2\omega(1)+2\gamma+4.
\end{equation}
Suppose that $0<t<1/(\sqrt{2}e^2)$ and choose
\begin{equation*}\label{eq:choices}
p=\left\lceil \log\left(\frac{1}{2t^2}\right)\right\rceil\ge 2 \qquad\mathrm{and} \qquad k=\left\lceil \left(\frac{2\omega(1)+2\gamma+4}{\omega(\sqrt{2} t)}\right)^{2\log\left(\frac{1}{2t^2}\right)}\right\rceil.
\end{equation*}
These choices  ensure that  $k^{\frac{1}{p}} \omega(\sqrt{2} t)\ge  2\omega(1)+2\gamma+4$ and $(2t^2)^{\frac{1}{p}}\ge 1/e$, so~\eqref{eq:omega constraint} implies that
$$
\omega(\sqrt{2} t)\ge \frac{(2t^2)^{\frac{1}{p}}}{2}\ge \frac{1}{2e}.
$$
Thus $\liminf_{s\to 0}\omega(s)>0$, a contradiction.\qed

\begin{remark}
{\em In order to obtain an example of separable Banach spaces $(X',\|\cdot\|_{X'})$ and $(Y',\|\cdot\|_{Y'})$ that satisfy the conclusion of Theorem~\ref{thm:no ext}, replace the $\ell_\infty$ products in~\eqref{eq:X,Y} by $c_0$ products, i.e., define
$$
X'\eqdef \Big(\bigoplus_{p=2}^\infty \ell_2\Big)_{c_0}\qquad \mathrm{and}\qquad Y'\eqdef \Big(\bigoplus_{p=2}^\infty  \ell_p\Big)_{c_0}.
$$
If the initial $1$-net $\mathscr{M}\subset \ell_2$ is chosen so that $0\in \mathscr{M}$, then the set $\mathscr{N}'\eqdef X'\cap \prod_{p=2}^\infty \mathscr{M}$ is a $1$-net in $X'$ and the above proof of Theorem~\ref{thm:no ext} goes through in this modified setting without any other change. }
\end{remark}

\subsection*{Acknowledgements} I am grateful to Bill Johnson for suggesting that I consider   Rosendal's question and several helpful discussions on these topics. I also thank Christian Rosendal for sharing his work~\cite{Ros15} and helpful comments, in particular noting that the proof of Theorem~\ref{thm:no ext} also yields~\eqref{eq:no bounded distance}.

\bibliographystyle{abbrv}
\bibliography{net-ext}

\end{document}